\documentclass{article}
\usepackage{amsmath, amsfonts, amssymb}
\usepackage[english]{babel}

\begin{document}
\begin{center}
{\Large \textbf{Boundary-value problems with non-local condition for degenerate parabolic equations}}

\medskip

J.M.Rassias$^1$ , E.T.Karimov $^2$

\medskip

$^1$ Pedagogical department, Mathematics and Informatics section, National and Capodistrian University of Athens, Athens, Greece\\
$^2$ Institute of Mathematics, National University of Uzbekistan named after Mirzo Ulugbek, Tashkent, Uzbekistan\\

\medskip
E-mails: {\tt jrassias@primedu.uoa.gr}\ (J.M.Rassias), \\
{\tt erkinjon@gmail.com}\ (E.T.Karimov)
\end{center}
\bigskip

\textbf{Abstract}

In this work we deal with degenerate parabolic equations with three lines of degeneration. Using "a-b-c" method we prove the uniqueness theorems defining conditions to parameters. We show nontrivial solutions for considered problems, when uniqueness conditions to parameters, participating in the equations are not fulfilled.

\medskip

\textbf{Keywords}: {Degenerate parabolic equation, boundary-value problems with non-local initial conditions, classical "a-b-c" method. eigenvalues and eigenfunctions.}\\

\textbf{MSC 2000}: {35K20, 35K65}

\medskip

\section{Introduction}

Degenerate partial differential equations have many applications in practical problems. Even finding some particular solutions for these kinds of equations is interesting for specialists in numerical analysis.

We note works by Friedman, Gorkov [1,2] on finding solutions for degenerate parabolic equations. Regarding the usage of nonlocal conditions used in the present work for parabolic equations we cite the work by Shopolov [4]. There are some papers [3,6-12] containing interesting results on investigation various type of degenerate partial differential equations. For detailed information on degenerate partial differential equations one can find in the monograph [5].

This work is continuation of the work [13] and only for convenience of the reader we did not omit some details of proofs.

\section{Boundary problem with nonlocal condition for parabolic equations with three lines of degeneration}

Let $\Omega $ be a simple-connected bounded domain in $R^3$ with boundaries $S_i$ ($i = \overline{1, 6}$). Here
$$
\begin{array}{l}
 S_1  = \left\{ {\left( {x,y,t} \right):\,t = 0,\,0 < x < 1,\,0 < y < 1} \right\},\\
 S_2  = \left\{ {\left( {x,y,t} \right):\,x = 1,\,0 < y < 1,\,0 < t < 1} \right\}, \\
 S_3  = \left\{ {\left( {x,y,t} \right):\,y = 0,\,0 < x < 1,\,0 < t < 1} \right\},\\
 S_4  = \left\{ {\left( {x,y,t} \right):\,x = 0,\,0 < y < 1,\,0 < t < 1} \right\}, \\
 S_5  = \left\{ {\left( {x,y,t} \right):\,y = 1,\,0 < x < 1,\,0 < t < 1} \right\},\\
 S_6  = \left\{ {\left( {x,y,t} \right):\,t = 1,\,0 < x < 1,\,0 < y < 1} \right\}. \\
 \end{array}
 $$

We consider a degenerate parabolic equation
\begin{equation}
x^n y^m u_t  = t^k y^m u_{xx}  + t^k x^n u_{yy}  - \lambda t^k x^n y^m u \label{eq.1}
\end{equation}
in the domain $\Omega $. Here $m > 0,\,n > 0,\,k>0,\,\lambda  = \lambda _1  + i\lambda _2 ,\,\,\lambda _1 ,\lambda _2  \in R$.

\textbf{The problem 1.} To find a function $u\left( {x,y,t} \right)$
 satisfying the following conditions:
\begin{enumerate}
\item $u\left( {x,y,t} \right) \in C\left( {\overline \Omega  } \right) \cap C_{x,y,t}^{2,2,1} \left( \Omega  \right)$;

\item $u\left( {x,y,t} \right)$ satisfies the equation (\ref{eq.1}) in $\Omega $;

\item $u\left( {x,y,t} \right)$ satisfies boundary conditions

\begin{equation}
u\left( {x,y,t} \right)\left| {_{S_2  \cup S_3  \cup S_4  \cup S_5 }  = 0} \right.;\label{eq.2}
\end{equation}

\item and a non-local condition
\begin{equation}
u\left( {x,y,0} \right) = \alpha u\left( {x,y,1} \right).\label{eq.3}
\end{equation}
\end{enumerate}
Here $\alpha  = \alpha _1  + i\alpha _2 ,\,\,\alpha _1 ,\alpha _2 $ are real numbers, moreover $\alpha _1^2  + \alpha _2^2  \ne 0$.

\textbf{Theorem 1.} If $\alpha _1^2  + \alpha _2^2  < 1,\,\,\lambda _1  \ge 0$ and exists a solution of the problem 1, then it is unique.

\verb"Proof":
Let us suppose that the problem 1 has two solutions $u_1 ,\,u_2 $. Denoting $u = u_1  - u_2 $
 we claim that $u \equiv 0$ in $\Omega $.

First we multiply equation (\ref{eq.1}) to the function $\overline u \left( {x,y,t} \right)$, which is complex conjugate function of $u\left( {x,y,t} \right)$. Then integrate it along the domain $\Omega _\varepsilon$ with boundaries
$$
\begin{array}{l}
 S_{1\varepsilon }  = \left\{ {\left( {x,y,t} \right):\,t = \varepsilon ,\,\varepsilon  < x < 1 - \varepsilon ,\,\varepsilon  < y < 1 - \varepsilon } \right\},\\
 S_{2\varepsilon }  = \left\{ {\left( {x,y,t} \right):\,x = 1 - \varepsilon ,\,\varepsilon  < y < 1 - \varepsilon ,\,\varepsilon  < t < 1 - \varepsilon } \right\}, \\
 S_{3\varepsilon }  = \left\{ {\left( {x,y,t} \right):\,y = \varepsilon ,\,\varepsilon  < x < 1 - \varepsilon ,\,\varepsilon  < t < 1 - \varepsilon } \right\},\\
 S_{4\varepsilon }  = \left\{ {\left( {x,y,t} \right):\,x = \varepsilon ,\,\varepsilon  < y < 1 - \varepsilon ,\,\varepsilon  < t < 1 - \varepsilon } \right\}, \\
 S_{5\varepsilon }  = \left\{ {\left( {x,y,t} \right):\,y = 1 - \varepsilon ,\,\varepsilon  < x < 1 - \varepsilon ,\,\varepsilon  < t < 1 - \varepsilon } \right\},\\
 S_{6\varepsilon }  = \left\{ {\left( {x,y,t} \right):\,t = 1 - \varepsilon ,\,\varepsilon  < x < 1 - \varepsilon ,\,\varepsilon  < y < 1 - \varepsilon } \right\}. \\
 \end{array}
$$
Then taking real part of the obtained equality and considering
$$
{\mathop{\rm Re}\nolimits} \left( {t^k y^m \overline u u_{xx} } \right) = {\mathop{\rm Re}\nolimits} \left( {t^k y^m \overline u u_x } \right)_x  - t^k y^m \left| {u_x } \right|^2 ,
$$
$$
{\mathop{\rm Re}\nolimits} \left( {t^k x^n \overline u u_{yy} } \right) = {\mathop{\rm Re}\nolimits} \left( {t^k x^n \overline u u_y } \right)_y  - t^k x^n \left| {u_y } \right|^2 ,
$$
$$
{\mathop{\rm Re}\nolimits} \left( {x^n y^m \overline u u_t } \right) = \left( {\frac{1}{2}x^n y^m \left| u \right|^2 } \right)_t ,
$$
after using Green's formula we pass to the limit at $\varepsilon  \to 0$. Then we get
$$
\begin{array}{l}
\iint\limits_{\partial\Omega}{\mathop{\rm Re}\nolimits} \left[ {t^k y^m \overline u u_x \cos \left( {\nu ,x} \right) + t^k x^n \overline u u_y \cos \left( {\nu ,y} \right) - \frac{1}{2}x^n y^m \left| u \right|^2 \cos \left( {\nu ,t} \right)} \right]d\tau  \\
  = \iiint\limits_{\Omega}\left( {t^k y^m \left| {u_x } \right|^2  + t^k x^n \left| {u_y } \right|^2  + \lambda _1 t^k x^n y^m \left| u \right|^2} \right)d\sigma,  \\
 \end{array}
 $$
where $\nu $ is outer normal.
Considering
$
{\mathop{\rm Re}\nolimits} \left[ {\overline u u_x } \right] = {\mathop{\rm Re}\nolimits} \left[ {u\overline u _x } \right],\,\,{\mathop{\rm Re}\nolimits} \left[ {\overline u u_y } \right] = {\mathop{\rm Re}\nolimits} \left[ {u\overline u _y } \right],
$
we obtain
$$
{\mathop{\rm Re}\nolimits} {\iint\limits_{S_1 } {\frac{1}{2}x^n y^m \left| u \right|^2 d\tau _1 } }  + {\iint\limits_{S_2 } {t^k y^m {\mathop{\rm Re}\nolimits} \left[ {u\overline u _x } \right]d\tau _2 } }  - {\iint\limits_{S_3 } {t^k x^n {\mathop{\rm Re}\nolimits} \left[ {u\overline u _y } \right]d\tau _3 } }
$$
\begin{equation}
- {\iint\limits_{S_4 } {t^k y^m {\mathop{\rm Re}\nolimits} \left[ {u\overline u _x } \right]d\tau _4 } }  + {\iint\limits_{S_5 } {t^k x^n {\mathop{\rm Re}\nolimits} \left[ {u\overline u _y } \right]d\tau _5 } }  - {\mathop{\rm Re}\nolimits} {\iint\limits_{S_6 } {\frac{1}{2}x^n y^m \left| u \right|^2 d\tau _6  = } }\label{eq.4}
\end{equation}
$$
 ={\iiint\limits_\Omega  {{\left( {t^k y^m \left| {u_x } \right|^2  + t^k x^n \left| {u_y } \right|^2  + \lambda _1 t^k x^n y^m \left| u \right|^2} \right)d\sigma } } }.
$$
From (\ref{eq.4}) and by using conditions (\ref{eq.2}), (\ref{eq.3}), we find
$$
\frac{1}{2}\left[ {1 - \left( {\alpha _1^2  + \alpha _2^2 } \right)} \right]\int\limits_0^1 {\int\limits_0^1 {x^n y^m \left| {u\left( {x,y,1} \right)} \right|^2 dxdy}  +}
$$
\begin{equation}
 +{\iiint\limits_\Omega  {{\left( {t^k y^m \left| {u_x } \right|^2  + t^k x^n \left| {u_y } \right|^2  + \lambda _1 t^k x^n y^m \left| u \right|^2} \right)d\sigma }  } }=0 .\label{eq.5}
\end{equation}
Setting $\alpha _1^2  + \alpha _2^2  < 1,\,\,\lambda _1  \ge 0$, from (\ref{eq.4}) we have $u\left( {x,y,t} \right) \equiv 0$
in $\overline \Omega  $.\\
Theorem 1 is proved.

We find below non-trivial solutions of the problem 1 at some values of parameter $\lambda $
 for which the uniqueness condition ${\mathop{\rm Re}\nolimits} \lambda  = \lambda _1  \ge 0$
 is not fulfilled.

We search the solution of Problem 1 as follows
\begin{equation}
u\left( {x,y,t} \right) = X\left( x \right) \cdot Y\left( y \right) \cdot T\left( t \right).\label{eq.6}
\end{equation}
After some evaluations we obtain the following eigenvalue problems:

\begin{equation}
\left\{
\begin{array}{l}
 X''\left( x \right) + \mu _1 x^n X\left( x \right) = 0, \\
 X\left( 0 \right) = 0,\,\,\,X\left( 1 \right) = 0;
\end{array}
\right.
\label{eq.7}
\end{equation}

\begin{equation}
\left\{
\begin{array}{l}
 Y''\left( y \right) + \mu _2 y^m Y\left( y \right) = 0, \\
 Y\left( 0 \right) = 0,\,\,\,Y\left( 1 \right) = 0;
\end{array}
\right.
\label{eq.8}
\end{equation}

\begin{equation}
\left\{ \begin{array}{l}
 T'\left( t \right) + \left( {\lambda  + \mu } \right)t^k T\left( t \right) = 0, \\
 T\left( 0 \right) = \alpha T\left( 1 \right). \\
 \end{array} \right.\label{eq.9}
\end{equation}
Here $\mu  = \mu _1  + \mu _2 $ is a Fourier constant.

Solving eigenvalue problems (\ref{eq.7}), (\ref{eq.8}) we find
\begin{equation}
\mu _{1l}  = \left( {\frac{{n + 2}}{2}\widetilde{\mu _{1l} }} \right)^2 ,\,\,\,\,\mu _{2p}  = \left( {\frac{{m + 2}}{2}\widetilde{\mu _{2p} }} \right)^2 ,\label{eq.10}
\end{equation}
\begin{equation}
X_l \left( x \right) = A_l \left( {\frac{2}{{n + 2}}} \right)^{\frac{1}{{n + 2}}} \mu _{1l}^{\frac{1}{{2\left( {n + 2} \right)}}} x^{\frac{1}{2}} J_{\frac{1}{{n + 2}}} \left( {\frac{{2\sqrt {\mu _{1l} } }}{{n + 2}}x^{\frac{{n + 2}}{2}} } \right),\label{eq.11}
\end{equation}
\begin{equation}
Y_p \left( y \right) = B_p \left( {\frac{2}{{m + 2}}} \right)^{\frac{1}{{m + 2}}} \mu _{2p}^{\frac{1}{{2\left( {m + 2} \right)}}} y^{\frac{1}{2}} J_{\frac{1}{{m + 2}}} \left( {\frac{{2\sqrt {\mu _{2p} } }}{{m + 2}}x^{\frac{{m + 2}}{2}} } \right),\label{eq.12}
\end{equation}
where $l,p = 1,2,...$, $\widetilde{\mu _{1l} }$ and $\widetilde{\mu _{2p} }$ are roots of equations $J_{\frac{1}{{n + 2}}} \left( x \right) = 0$
 and $J_{\frac{1}{{m + 2}}} \left( y \right) = 0$, respectively, $J_p (\cdot)$ is the first kind Bessel function of p-th order.

The eigenvalue problem (\ref{eq.9}) has non-trivial solution only when
$$
\left\{ \begin{array}{l}
 \frac{\alpha _1}{\alpha_1^2+\alpha_2^2}  = e^{-\frac{\lambda _1  + \mu _{lp} }{k+1}} \cos \frac{\lambda_2}{k+1},  \\
 \frac{\alpha _2}{\alpha_1^2+\alpha_2^2}  = e^{-\frac{\lambda _1  + \mu _{lp} }{k+1}} \sin \frac{\lambda_2}{k+1} . \\
 \end{array} \right.
$$

Here $\lambda  = \lambda _1  + i\lambda _2 ,\,\,\,\,\alpha  = \alpha _1  + i\alpha _2 ,\,\,\,\mu _{lp}  = \mu _{1l}  + \mu _{2p} $. After elementary calculations, we get
\begin{equation}
\lambda _1  =  - \mu _{lp}  + \frac{k+1}{2}\ln \left(\alpha _1^2  + \alpha _2^2 \right) ,\,\,\,\,\lambda _2  = (k+1)\left [\arctan \frac{{\alpha _2 }}{{\alpha _1 }} + s\pi\right ] ,\,\,\,s \in Z^+.\label{eq.13}
\end{equation}
Corresponding eigenfunctions have the form
\begin{equation}
T_s \left( t \right) = C_s e^{\left[ {- \ln \sqrt {\alpha _1^2  + \alpha _2^2 }  - i\left( {\arctan \frac{{\alpha _2 }}{{\alpha _1 }} + s\pi } \right)} \right]t^{k+1}} .\label{eq.14}
\end{equation}

Considering (\ref{eq.6}), (\ref{eq.11}), (\ref{eq.12}) and (\ref{eq.14}) we can write non-trivial solutions of the problem 1 in the following form:
$$
u_{lps} \left( {x,y,t} \right) = D_{lps} \left( {\frac{2}{{n + 2}}} \right)^{\frac{1}{{n + 2}}} \left( {\frac{2}{{m + 2}}} \right)^{\frac{1}{{m + 2}}} \mu _{1l}^{\frac{1}{{2\left( {n + 2} \right)}}} \mu _{2p}^{\frac{1}{{2\left( {m + 2} \right)}}} \sqrt {xy} $$
$$
 \times J_{\frac{1}{{n + 2}}} \left( {\frac{{2\sqrt {\mu _{1l} } }}{{n + 2}}x^{\frac{{n + 2}}{2}} } \right)\cdot J_{\frac{1}{{m + 2}}} \left( {\frac{{2\sqrt {\mu _{2p} } }}{{m + 2}}y^{\frac{{m + 2}}{2}} } \right)
 e^{\left[ {- \ln \sqrt {\alpha _1^2  + \alpha _2^2 }  - i\left( {\arctan \frac{{\alpha _2 }}{{\alpha _1 }} + s\pi } \right)} \right]t^{k+1}},
 $$ 									
where $D_{lps}  = A_l  \cdot B_p  \cdot C_s $ are constants.

\textbf{Remark 1.} One can easily see that $\lambda _1  < 0$
 in (\ref{eq.13}), which contradicts to condition $Re{\lambda}  = \lambda _1  \ge 0$
 of the Theorem 1.

\textbf{Remark 2.} The following problems can be studied by similar way.
Instead of condition (\ref{eq.2}) we put zero conditions on surfaces as follows:

\begin{tabular}{|c|c|c|c|c|c|c|c|c|}
  \hline
  Problem's name & P$_2$ & P$_3$ & P$_4$ & P$_5$ & P$_6$ & P$_7$ & P$_8$ & P$_9$ \\
  \hline
  S$_2$ & $u_x$ & $u$ & $u$ & $u_x$ & $u$ & $u$ & $u_x$ & $u$ \\
  S$_3$ & $u_y$ & $u$ & $u_y$ & $u$ & $u_y$ & $u$ & $u$ & $u$ \\
  S$_4$ & $u$ & $u_x$ & $u$ & $u_x$ & $u$ & $u_x$ & $u$ & $u$ \\
  S$_5$ & $u$ & $u_y$ & $u_y$ & $u$ & $u$ & $u$ & $u$ & $u_y$ \\
  \hline
\end{tabular}

\section{Boundary problem with two nonlocal conditions}
We consider equation
\begin{equation}
t^ky^mU_{xx}-t^kx^nU_y-x^ny^mU_t-\Lambda x^ny^mt^kU=0 \label{eq.15}
\end{equation}
in the domain $\Omega$. Here $m,n,k>0,\,\Lambda=\Lambda_{11}+\Lambda_{21}+i\left(\Lambda_{12}+\Lambda_{22}\right),\,\Lambda_{jj}\in \textbf{R}$.

\textbf{The Problem A. }\emph{To find a function $U(x,y,t)$ from the class of
$$
W=\left\{U(x,y,t): U\in C\left(\overline{\Omega}\right)\cap C_{x,y,t}^{2,1,1}\left(\Omega\right)\right\}
$$
satisfying equation (\ref{eq.15}) in $\Omega$
 and boundary condition}
\begin{equation}
U\left( {x,y,t} \right)\left| {_{S_2  \cup S_3 }  = 0} \right.,
\end{equation}
\emph{nonlocal conditions}
\begin{equation}
U\left( {x,0,t} \right) = \beta U\left( {x,1,t} \right),
\end{equation}
\begin{equation}
U\left( {x,y,0} \right) = \gamma U\left( {x,y,1} \right).
\end{equation}

Here $\gamma=\gamma_1+i\gamma_2,\,\beta=\beta_1+i\beta_2$.

\textbf{Theorem 2.}
If $\beta_1^2+\beta_2^2<1,\,\,\gamma_1^2  + \gamma_2^2  < 1,\,\,\Lambda_{11}+\Lambda_{21}\ge 0$ and exists a solution of the problem A, then it is unique.

Proof of this theorem will be done similarly as in the Theorem 1.
Nontrivial solutions of this problem can be written as follows:
$$
U_{ls} \left( {x,y,t} \right) = E_{ls} \left( {\frac{2}{{n + 2}}} \right)^{\frac{1}{{n + 2}}} \mu _{1l}^{\frac{1}{{2\left( {n + 2} \right)}}} \sqrt {x} J_{\frac{1}{{n + 2}}} \left( {\frac{{2\sqrt {\mu _{1l} } }}{{n + 2}}x^{\frac{{n + 2}}{2}} } \right)
$$
$$
 \times e^{\left[ {- \ln \sqrt {\beta_1^2  + \beta_2^2 }  - i\left( {\arctan \frac{{\beta_2 }}{{\beta_1 }} + s\pi } \right)} \right]y^{m+1}}
 e^{\left[ {- \ln \sqrt {\gamma_1^2  + \gamma_2^2 }  - i\left( {\arctan \frac{{\gamma_2 }}{{\gamma_1 }} + s\pi } \right)} \right]t^{k+1}},
 $$ 									
where $E_{ls}$ are constants.

We note that above given nontrivial solutions exist only when
$$
\left\{ \begin{array}{l}
 \frac{\beta_1}{\beta_1^2+\beta_2^2}  = e^{-\frac{\Lambda_{11}  + \mu_{1l}}{m+1}} \cos \frac{\Lambda_{12}}{m+1},  \\
 \frac{\beta_2}{\beta_1^2+\beta_2^2}  = e^{-\frac{\Lambda_{11}  + \mu _{1l} }{m+1}} \sin \frac{\Lambda_{12}}{m+1},\\
 \frac{\gamma_1}{\gamma_1^2+\gamma_2^2}  = e^{-\frac{\Lambda_{21}  + \mu _{1l} }{k+1}} \cos \frac{\Lambda_{22}}{k+1},  \\
 \frac{\gamma_2}{\gamma_1^2+\gamma_2^2}  = e^{-\frac{\Lambda_{21}  + \mu _{1l} }{k+1}} \sin \frac{\lambda_{22}}{k+1} . \\
 \end{array} \right.
$$
Here
$$
\Lambda_{11}=-\mu_{1l}+(m+1)ln\sqrt{\beta_1^2+\beta_2^2},\,\Lambda_{12}=(m+1)\arctan{\frac{\beta_2}{\beta_1}},
$$
$$
\Lambda_{21}=-\mu_{1l}+(k+1)ln\sqrt{\gamma_1^2+\gamma_2^2},\,\Lambda_{22}=(k+1)\arctan{\frac{\gamma_2}{\gamma_1}}.
$$

\end{document}